\documentclass[12pt]{amsart}
\usepackage[margin=1in]{geometry}
\usepackage[english]{babel}
\usepackage[utf8]{inputenc}
\usepackage{subcaption}
\usepackage{amsmath}
\usepackage{amssymb}
\usepackage{listings}
\usepackage{url}
\usepackage{bm}
\usepackage{graphicx}
\usepackage{psfrag}
\usepackage[dvipsnames]{xcolor}
\usepackage{color}
\usepackage{hyperref}
\hypersetup{colorlinks=true, linkcolor=blue, citecolor=magenta, urlcolor=black}

\definecolor{codegreen}{rgb}{0,0.6,0}\definecolor{codegray}{rgb}{0.5,0.5,0.5}
\definecolor{codepurple}{rgb}{0.58,0,0.82}
\definecolor{backcolour}{rgb}{0.95,0.95,0.92}

\lstdefinestyle{custom}{
    backgroundcolor=\color{backcolour},   
    commentstyle=\color{codegreen},
    keywordstyle=\color{red},
    stringstyle=\color{blue},
    basicstyle=\ttfamily\footnotesize,
    breakatwhitespace=false,         
    breaklines=true,                 
    captionpos=b,                    
    keepspaces=true,                 
    numbers=left,                    
    numbersep=5pt,                  
    showspaces=false,                
    showstringspaces=false,
    showtabs=false,                  
    tabsize=2
}

\newcommand{\RR}{\mathbb{R}}
\newcommand{\ZZ}{\mathbb{Z}}
\newcommand{\bfp}{\mathbf{p}}
\newcommand{\bfu}{\mathbf{u}}
\newcommand{\bfv}{\mathbf{v}}

\theoremstyle{plain}
\newtheorem{theorem}{Theorem}[section]
\newtheorem{lemma}[theorem]{Lemma}

\theoremstyle{definition}

\theoremstyle{remark}


\begin{document}
\title{On lattice triangles satisfying $\boldsymbol{B(T)=3}$ with collinear interior lattice points}
\author{Eddy Li}
\address{Department of Mathematics, Stanford Online High School, 415 Broadway Academy Hall, Redwood City, CA 94063 USA}
\email{eddyli@ohs.stanford.edu}
\author{Dana Paquin}
\address{Department of Mathematics, California Polytechnic State University, 1 Grand Ave, San Luis Obispo, CA 93407 USA}
\email{dpaquin@calpoly.edu}
\date{14 January 2025}
\bibliographystyle{alpha}
\renewcommand\qedsymbol{$\blacksquare$}
\begin{abstract}
A lattice point in $\RR^2$ is a point $(x,y)$ with $x,y\in\ZZ$, and a lattice triangle is a triangle whose three vertices are all lattice points. We investigate the integers $k$ with the property that if $T$ is a lattice triangle with $3$ boundary points and $k$ points in the interior, then all $k$ boundary points must be collinear.
\end{abstract}
\maketitle
%\tableofcontents
\section{Introduction}\label{sec:Intr}

For any positive integer $n$, a \textit{lattice point} is an element of $\mathbb{Z}^n$ embedded in the space $\mathbb{R}^n$. A \textit{lattice polytope} is a convex hull of finitely many lattice points, and a \textit{lattice polygon} is a lattice polytope when $n=2$. A \textit{lattice line segment} is a line segment whose endpoints are lattice points. For lattice points $V_1$ and $V_2$, let $L(V_1V_2)$ be the number of lattice points on the segment $V_1V_2$, excluding the endpoints. Additionally, for any lattice polytope $P$, let $B(P)$, $I(P)$, and $A(P)$ denote the number of boundary lattice points of $P$, the number of interior lattice points of $P$, and the area of $P$, respectively. Pick's Theorem \cite{pick} states that if $P$ is any lattice polygon, convex or otherwise, then $$A(P)=\frac{B(P)}{2}+I(P)-1.$$ A classic question in geometric number theory, or lattice point geometry, regards the generalization of Pick's Theorem to three dimensions \cite{si,yamamoto} and beyond \cite{hofscheier2,hofscheier3,hofscheier}.

A closely related question in geometric number theory, and similarly classic, considers the range of possible constructions of lattice polytopes with various configurations of interior and boundary lattice points. For instance, Hensley \cite{hensley} showed that the value of $B(P)$ for some lattice polytope $P$ satisfying $I(P)>0$ is bounded above by a function of $n$ and $I(P)$. Averkov, Kr\"umpelmann, and Nill \cite{averkov} extended Hensley's results by establishing the maximal volume of a lattice simplex $P$ satisfying $I(P)=1$, while Balletti and Kasprzyk \cite{balletti} obtained empirical evidence in support for an analogous volume bound when $I(P)=2$.

The analogous problem concerning the configurations of lattice points for lattice polygons in two dimensions is much more tractable. For example, if $T$ is a lattice triangle satisfying $I(T)=0$, we can show by direct construction that any number of boundary lattice points is possible. In particular, the lattice triangle with vertices $(0,0)$, $(1,0)$, and $(0,k-2)$ has $I=0$ interior lattice points and $B=k$ boundary lattice points.  On the other hand, using number theoretic properties of the greatest common divisor, properties of divisibility, and Pick's Theorem, we can show that if $T$ is a lattice triangle such that $I(T)=1$, then the only possible values of $B(T)$ are $3$,~$4$,~$6$,~$8$,~or~$9$ \cite{weaver}, and all of these values can be obtained via direct construction. 

In \cite{wei}, Wei and Ding showed that any convex lattice polygon with $11$ vertices must contain at least $15$ interior lattice points, and in \cite{wei2}, the same authors classify, up to lattice equivalence, all lattice polygons with exactly $2$ interior lattice points. Arkinstall \cite{arkinstall} proved that there is only one convex lattice hexagon with exactly one interior lattice point, again up to lattice equivalence, and in \cite{rabinowitz,rabinowitz2}, Rabinowitz extended the result by finding all convex lattice polygons with at most one interior lattice point.

In the $23$rd Bay Area Mathematical Olympiad \cite{bamo}, the fifth problem was to prove that if $T$ is a lattice triangle with $B(T)=3$ and $I(T)=4$, then the $4$ interior lattice points of $T$ must all be collinear. In this paper, we extend that result to classify all of those positive integers $k$ with the property that all lattice triangles $T$ satisfying $B(T)=3$ and $I(T)=k$ must have its $k$ interior lattice points be collinear. For simplicity of notation, we will call an integer $k$ with this property a \textit{$2$-collinear} integer, following \cite{reznick}, where the $2$ here represents the dimension of the simplex -- in our case, a triangle in $\mathbb{Z}^2$.  In particular, we will show that the only $2$-collinear integers are $1$, $2$, $4$, and $7$.
 
This paper is organized in the following way. In Section 2, we will establish the relationship between lattice triangles and unimodular shear transformations. In Section 3, we will prove that $1$, $2$, $4$, and $7$ are all $2$-collinear, integers spending most of our time on the case $k=7$, and also determine criteria for integers that are not $2$-collinear. Finally, in Section 4, we use the Schemmel totient function, a number-theoretic variation of the ordinary Euler totient function, to show that the {\em only} $2$-collinear integers are $1$, $2$, $4$, and $7$. Section 5 concludes.

\section{Integer Unimodular Shear Transformations}

A linear transformation $f:\RR^2\to\RR^2$ is a \textit{unimodular shear} if the standard matrix representation $M$ of $f$ satisfies $$|\!\det M|=1.$$ An \textit{integer unimodular shear} is a unimodular shear with the additional property that the standard matrix representation $M$ of $f$ has integer entries.   If $T$ is a triangle with vertices $V_1$, $V_2$, and $V_3$, we let $f(T)$ denote the triangle with vertices $f(V_1)$, $f(V_2)$, and $f(V_3)$, and if $V_1V_2$ is a line segment with endpoints $V_1$ and $V_2$, we let $f(V_1V_2)$ denote the line segment with endpoints $f(V_1)$ and $f(V_2)$. 

To begin, we establish some basic properties of integer unimodular shears with the following two lemmas.

\begin{lemma}\label{shearinvert}  If $f$ is an integer unimodular shear transformation, then $f$ is invertible and the inverse $g$ of $f$ is an integer unimodular shear.\end{lemma}
\begin{proof}Let $f$ be an integer unimodular shear transformation, and let
\begin{equation*}
M=\begin{bmatrix}a&b\\c&d\end{bmatrix}
\end{equation*}
be the standard matrix representation of $f$.  Since $f$ is an integer unimodular shear, $a,b,c,d\in\mathbb{Z}$ and $\text{det} M=1.$ Consider the linear transformation $g$ with standard matrix representation \begin{equation*}
N=\begin{bmatrix}d&-b\\-c&a\end{bmatrix}.
\end{equation*}
Since $f$ is unimodular, $$|ad-bc|=|\!\det M|=1,$$ so $$MN=\begin{bmatrix}a&b\\c&d\end{bmatrix}\begin{bmatrix}d&-b\\-c&a\end{bmatrix}=\begin{bmatrix}ad-bc&0\\0&ad-bc\end{bmatrix}=\begin{bmatrix}1&0\\0&1\end{bmatrix}=I,$$ where $I$ is the identity matrix. Similarly, $NM=I$. Thus, $g$ is the inverse transformation of $f$.  Further, we have $$|\!\det N|=|da-(-b)(-c)|=|ad-bc|=1,$$ and the entries of $N$ are integers by construction. Since the standard matrix of $g$ is $N$, it follows that $g$ is an integer unimodular shear.  We conclude that $f$ is invertible over the set of integer unimodular shears. \end{proof}

\begin{lemma}\label{shearline}Collinearity is preserved under the action of an integer unimodular shear. Additionally, if $f$ is an integer unimodular shear and if $V_1V_2$ is a lattice line segment, then $$L(V_1V_2)=L(f(V_1V_2).$$\end{lemma}
\begin{proof}We first prove that $f$ preserves collinearity. In particular, note that any line $\ell$ can be parameterized as the set of points corresponding to vectors of the form $\bfu+\bfv t$, where $\bfu$ and $\bfv$ are fixed vectors in $\RR^2$ and $t$ varies over $\RR$.

Now fix some $t_0$, and consider the point  $\bfp$ on $\ell$ given by $\bfp=\bfu+\bfv t_0$. Let $f$ be an integer unimodular shear with standard matrix representation $M$.  Then $f$ maps $\bfp=\bfu+\bfv t_0$ to the point $$f(\bfp)=M\bfp=M(\bfu+\bfv t_0)=(M\bfu)+(M\bfv)t_0,$$ which must lie on the line $\ell'$ parameterized by the vector expression $(M\bfu)+(M\bfv)t$. Since $\bfu$ and $\bfv$ are fixed, it thus follows that any points on $\ell$ must be mapped onto the line $\ell'$ by $f$. Hence $f$ preserves collinearity.

Next, we show that $f$ preserves the value of $L$ when transforming some lattice segment $V_1V_2$. Using Lemma~\ref{shearinvert}, let $g$ and $N$ be the shear transformation that inverts $f$ and its corresponding standard matrix, respectively.

By definition of integer unimodular shear, we see that each lattice point on $V_1V_2$ is mapped to another lattice point on its image, since $M$ has integer entries. Similarly, the converse must also be true since $N$ has integer entries. Thus $M$ and $N$ induce a bijection between the lattice points on $V_1V_2$ and $f(V_1V_2)$, yielding the equality $L(V_1V_2)=L(f(V_1V_2)).$ \end{proof}
Lemma~\ref{shearline} gives rise to the following result.

\begin{lemma}\label{sheartriangle}Let $f$ be an integer unimiodular shear, and let $T$ be a lattice triangle. Then $f$ preserves the values of $B$, $I$, and $A$ between the original triangle $T$ and its image $f(T)$.\end{lemma}
\begin{proof}First we prove that $B(T)=B(f(T)).$ Let the vertices of $T$ be $V_1$, $V_2$, and $V_3$ Then $$B(T)=L(V_1V_2)+L(V_2V_3)+L(V_3V_1)+3.$$ By Lemma~\ref{shearline}, 
\begin{align*}
B(T)&=L(V_1V_2)+L(V_2V_3)+L(V_3V_1)+3\\
&=L(f(V_1V_2))+L(f(V_2V_3))+L(f(V_3V_1))+3=B(f(T)).
\end{align*}

Next we prove that $A(T)=A(f(T)).$ Here, since $|\!\det(M)|=1$ by definition, it follows that $A(f(T))=A(T)|\!\det(M)|=A(T)$.

Finally, by Pick's Theorem, we have $I(T)=A(T)-\frac{B(T)}2+1,$ so that $$I(I)=A(T)-\frac{B(T)}2+1=A(f(T))-\frac{B(f(T))}2+1=I(f(T)).$$ Hence, we conclude that $f$ preserves the values of $B$, $I$, and $A$ for $T$ and $f(T).$
\end{proof}
We will use Lemma~\ref{sheartriangle} to obtain a much simpler but equivalent method of determining whether an integer $k$ is $2$-collinear.

\begin{theorem}\label{goodk}An integer $k$ is $2$-collinear if and only if for all triangles $T$ such that two of its vertices are $(0,0)$ and $(1,0)$, and $(B(T),I(T))=(3,k)$, the $k$ interior points are collinear.\end{theorem}
\begin{proof}Since integer unmiodular shears preserve the values of $B$, $I$, $A$, as well as the property of collinearity, due to Lemmas~\ref{shearline} and~\ref{sheartriangle}, it suffices to prove that we can map any triangle $T$ with $B(T)=3$ and $I(T)=k$ another triangle $T'$ with two of its vertices being $(0,0)$ and $(1,0)$ using only integer unimodular shears and lattice translations.

First we translate $T$ so that one of its vertices land on the origin; let $V_1=(x_1,y_1)$ and $V_2=(x_2,y_2)$ be the other two vertices. It then suffices to shear $(x_1,y_1)$ onto $(1,0).$

Note that $B(T)=3$ implies that there can be no lattice points on the sides of $T$; otherwise we would have $B(T)\ge 4.$ Now suppose that $\gcd(x_1,y_1)=d>2$. Then $(\frac{x_1}{d},\frac{y_1}{d})$ would be on $V_1O$, so that $L(V_1O)\ge1,$ a contradiction, forcing $\gcd(x_1,y_1)$ to be $1.$ 

Thus, Bezout's Lemma implies that $cx_1+dy_1=\gcd(c,d)=1$ for some choice of $c$ and $d$. Now consider the transformation with standard matrix $$M=\begin{bmatrix}c&d\\y_1&-x_1\end{bmatrix},$$ which is an integral unimodular shear as $$\det(M)=c(-x_1)-d(y_1)=-1.$$ Applying $M$ then maps $(x_1,y_1)$ onto the lattice point $(cx_1+dy_1,y_1x_1-x_1y_1)=(1,0)$, as desired.\end{proof}

\section{2-Collinear Integers}

In this section, we prove that $1$, $2$, $4$, and $7$ are all $2$-collinear integers and also establish criteria for integers that are not $2$-collinear.
\begin{theorem}\label{124}The integers $1$, $2$, and $4$ are all $2$-collinear.\end{theorem}
\begin{proof} Clearly $1$ and $2$ are both $2$-collinear integers. The statement that $4$ is a $2$-collinear integer was the fifth problem of the 2023 BAMO-12, with a proof given in \cite{bamo}.\end{proof}

Next, we will prove the following result.
\begin{lemma}\label{gcd}Lattice points $V_1=(x_1,y_1)$ and $V_2=(x_2,y_2)$ satisfy $L(V_1V_2)=0$ if and only if $\gcd(x_2-x_1,y_2-y_1)=1.$\end{lemma}
\begin{proof}First we prove the if direction. In particular, suppose that $\gcd(x_2-x_1,y_2-y_1)=1$. Then any point $V_3=(x_3,y_3)$ on the line segment strictly between $V_1$ and $V_2$ must satisfy $$\frac{x_3-x_1}{y_3-y_1}=\frac{x_2-x_1}{y_2-y_1},$$ with $|x_3-x_1|<|x_2-x_1|$ and $|y_3-y_1|<|y_2-y_1|$.

Let $$r=\frac{x_2-x_1}{x_3-x_1}=\frac{y_2-y_1}{y_3-y_1}=\frac{p}q,$$ with $\gcd(p,q)=1$, so that $r>1$ and thus $p>1$. In particular, it follows that $$\gcd(x_2-x_1,y_2-y_1)\ge p>1,$$ a contradiction. Thus $$L(V_1V_2)=0.$$

Now we prove the converse of the statement. Suppose that $\gcd(x_2-x_1,y_2-y_1)=d>1,$ and consider the point $$V_3=\left(x_1+\frac{x_2-x_1}d,y_1+\frac{y_2-y_1}d\right)=\left(\frac{(d-1)x_1+x_2}d,\frac{(d-1)y_1+y_2}d\right).$$
Clearly, $V_3$ is a lattice point by construction.

Now we prove that $V_3$ is on the line segment $V_1V_2$. Let $\bfv_i\in\RR^2$ be the vector with entries equal to that of $V_i.$ Since our segment is parallel to the vector $\bfv_2-\bfv_1$ going from $V_1$ to $V_2,$ we can thus parameterize it as being equivalent to the set of points represented by the equation $$\bfv_1+t(\bfv_2-\bfv_1)=(1-t)\bfv_1+t\bfv_2.$$ Setting $t=\frac{1}d$ then yields the position of $V_3$, which implies that $V_3$ is strictly between $V_1$ and $V_2$ as $0<t<1$, giving $L(V_1V_2)\ge1$, a contradiction. Thus such a $V_3$ cannot exist, so $d=1.$
\end{proof}

\begin{lemma}\label{shifts}Suppose that the triangle $T$ with vertices $O=(0,0)$, $W=(1,0)$, and $V=(a,b)$ has $I(T)$ collinear interior points. Then the triangle $T'$ with vertices at the lattice points $O$, $W$, and $V'=(a+tb,b)$ for $t\in\ZZ$ also have $I(T)$ collinear interior points for any integer $t.$\end{lemma}

\begin{proof}Consider the transformation $f$ with standard matrix $M$ given by $$M=\begin{bmatrix}1&t\\0&1\end{bmatrix}.$$ One can check that $f$ is a shear with the properties that $f(O)=O$, $f(W)=W$, and $f(V)=V'$, so that $f(T)=T'$. Then Lemma~\ref{sheartriangle} implies that $I(T)=I(T')$, as desired.\end{proof}

\begin{theorem}\label{7}The integer $7$ is $2$-collinear.\end{theorem}
\begin{proof}Using Theorem~\ref{goodk}, we only need to show that if $T$ is a triangle with vertices $O=(0,0)$, $W=(1,0)$, and $V=(a,b)$ such that $B(T)=3$ and $I(T)=7,$ then these seven interior points must all be collinear.

First, by Pick's Theorem, we have that $$A(T)=\frac{B(T)}2+I(T)-1=\frac32+7-1=\frac{15}2.$$ Also, the triangle area formula implies that $T$ has base $1$ and height $b$, so that $$\frac{b}{2}=\frac{15}2,$$ implying that $b=15.$ The fact that $B(P)=3$ implies that $$L(VO)=L(VW)=L(WO)=0,$$ so Lemma~\ref{gcd} implies that $$\gcd(a,b)=\gcd(a-1,b)=1;$$ equivalently, $$\gcd(a,15)=\gcd(a-1,15)=1.$$ Thus $$a\in \{2,8,14\},$$ and we verify directly that each of the triangles with vertices $$(0,0),~(1,0),~\text{ and }(a,15),$$ where $a\in\{2,8,14\}$ contain exactly $7$ interior lattice points, and these interior lattice points are all collinear.\end{proof}

Theorems~\ref{124} and~\ref{7} imply that $1$, $2$, $4$, and $7$ are all $2$-collinear integers. In particular, note that the argument from Theorem~\ref{7} can be generalized to show that $4$ is a $2$-collinear integer as well. We now have the following.
\begin{theorem}\label{1247}The integers $\{1,2,4,7\}$ are all $2$-collinear integers.\end{theorem}

We have thus proven that $1$, $2$, $4$, and $7$ are all $2$-collinear integers, though we have not determined whether or not these are the {\em only} $2$-collinear integers. In the following, we establish a sufficient condition for an integer $k$ to {\em not} satisfy the $2$-collinearity property.

\begin{theorem}\label{badk}Fix some integer $k.$ If there exists some positive integer $a$ such that $3\le a\le k$ and $\gcd(a,2k+1)=\gcd(a-1,2k+1)=1,$ then $k$ is not $2$-collinear.\end{theorem}

\begin{proof}Suppose that we have the above setup. Consider the triangle $T$ with vertices at the points $O=(0,0)$, $W=(1,0)$, and $V=(a,2k+1)$. By Lemma~\ref{gcd}, $$\gcd(a,2k+1)=\gcd(a-1,2k+1)=1$$ implies that $$L(VO)=L(VW)=L(WO)=0.$$ It follows that $B(T)=3.$

The triangle area formula implies that $A(T)=\frac{2k+1}2$, so Pick's Theorem implies that $$I(T)=A(T)-\frac{B(T)}2+1=\frac{2k+1}2-\frac32+1=k.$$ To prove that $k$ is not $2$-collinear, it suffices to prove that the $k$ interior points in $T$ are not all collinear.

Consider the intersection of $T$ with the line $x=1.$ Clearly one intersection point is $(1,0).$ The other point is given by $(1,\frac{2k+1}a).$ It follows that the line $x=1$ intersects $T$ at $p=\lfloor\frac{2k+1}a\rfloor$ lattice points; because $3\le a\le k,$ it follows that $\frac{2k+1}a$ is between $\frac{2k+1}k\geq 2$ and $\frac{2k+1}3<k$, so $$2\le p\le k-1.$$ Suppose, by way of contradiction, that the $k$ interior points in $T$ were all collinear. Since $p\ge 2,$ it follows that all $k$ points must lie on the line $x=1$; however, $p\le k-1$, a contradiction. Hence, not all of these points are collinear, so $k$ is not $2$-collinear.\end{proof}

\section{The Schemmel Totient Function}
In this section, we introduce the generalized totient functions and develop some essential features of this class of functions. We will then use the Schemmel totient function \cite{schemmel} to prove that $1$, $2$, $4$, and $7$ are the only $2$-collinear integers.

Let $k$ be an integer with $r$ distinct prime factors, so that $$k=p_1^{c_1}p_2^{c_2}\cdots p_r^{c_r}$$ for primes $p_i$ and positive integers $c_i.$ Then, due to \cite{alder}, define the \textit{generalized totient function} $\phi(k,m)$ to equal $$\phi(k,m)=k\prod_{i=1}^r\left(1-\frac{\epsilon(p_i,m)}{p_i}\right),$$ where $\epsilon(p_i,m)$ equals $1$ if $p_i\mid m$ and $2$ otherwise. It is easy to see that $\epsilon(p_i,0)=1$ for all $p_i$, implying that $\phi(k,0)=\phi(k)$, the familiar Euler totient function.

Now, suppose that we had set $m=1$ instead. We would then have $p_i\nmid 1$ and equivalently $\epsilon(p_i,1)=2$ for every $p_i$, so our function $\phi(k,m)$ would become $$\phi(k,1)=k\prod_{i=1}^r\left(1-\frac{2}{p_i}\right)=\left(\prod_{i=1}^r(p_i-2)\right)\left(\prod_{i=1}^rp_i^{c_i-1}\right).$$ This is known as the \textit{Schemmel totient function}, so named after its discoverer in \cite{schemmel}.

\begin{lemma}\label{multgt}If $a$ and $b$ are relatively prime positive integers, then $$\phi(a,m)\phi(b,m)=\phi(ab,m)$$ for all positive integers $m.$\end{lemma}
\begin{proof}The multiplicativity of generalized totient functions was shown in \cite{alder}.\end{proof}

We can use Lemma~\ref{multgt} to establish a loose bound on the Schemmel totient $\phi(k,1)$, which we will utilize to show that $1$, $2$, $4$, and $7$ are the only $2$-collinear integers.

\begin{lemma}\label{stg5}The Schemmel totient $\phi(k,1)$ is less than $5$ if and only if $k\in\{1,3,5,9,15\}$ or $2\mid k.$\end{lemma}
\begin{proof}First note that $\phi(1,1)=\phi(3,1)=1$, $\phi(5,1)=3$, $\phi(9,1)=3$, and $\phi(15,1)=3.$ Also, if $k$ is even, we let $p_1=2$ without loss of generality, so that
$$\phi(k,1)=\left(\prod_{i=1}^r(p_i-2)\right)\left(\prod_{i=1}^rp_i^{c_i-1}\right)=(2-2)\left(\prod_{i=2}^r(p_i-2)\right)\left(\prod_{i=1}^rp_i^{c_i-1}\right)=0.$$
It thus follows that if $k$ is even or $k\in\{1,3,5,9,15\},$ then $\phi(k,1)<5.$

Now suppose that $\phi(k,1)=0.$ Recalling the above summation definition of $\phi(k,1)$, it follows that the only possible way for this to happen is for some $p_i-2$ to equal $0$, so that $p_i=2$; in other words, $\phi(k,1)=0$ iff $2\mid k.$

Now, suppose that $k=p_1^{c_1}p_2^{c_2}\cdots p_r^{c_r}$ is odd and $0<\phi(k,1)<5.$ Assume that there exists some $p_j\ge 7,$ so it follows that $$\phi(k,1)=\left(\prod_{i=1}^r(p_i-2)\right)\left(\prod_{i=1}^rp_i^{c_i-1}\right)\ge\left(\prod_{i=1}^r(p_i-2)\right)\ge p_j-2\ge5,$$
a contradiction.

It thus follows that the maximum prime factor of $k$ is $5$; since $k$ is odd, the primes $3$ and $5$ are the only possible prime divisors of $k$, so $r$ is either $1$ or $2.$ If $r=2,$ then we let $p_1=3$ and $p_2=5$ without losing generality, so that $k=3^{c_1}5^{c_2}$ and $$\phi(k,1)=\left(1-\frac23\right)\left(1-\frac25\right)k=\frac{k}{5}.$$
If $r=2$ and $\phi(k,1)<5,$ then we are forced to take $c_1=c_2=1$, so that $k=15.$

Now suppose that $r=1,$ so that $\phi(k,1)=(p_1-2)p_1^{c_1-1}$ and $p_1\in\{3,5\}.$ If $p_1=3,$ then $\phi(k,1)=3^{c_1-1}<5,$ forcing $c_1$ to be either $1$ or $2$, respectively yielding the values $k=3$ and $k=9$. Otherwise, if $p_1=5,$ then $\phi(k,1)=3\cdot 5^{c_1-1}<5,$ so we must have $c_1=1$ and therefore $k=5.$ It follows that the only possible values of $k$ for this case are $3$, $5$, and $9$.

Since the above are the only possible cases for $r$, we have found all possible values of $k$, and the desired results immediately follows.
\end{proof}
Note that $3=2\cdot1+1$, $5=2\cdot2+1$, $9=2\cdot4+1$, and $15=2\cdot7+1$, so that the image of $\{1,2,4,7\}$ under the function $f(k)=2k+1$ is exactly $\{3,5,9,15\}.$ This is clearly not a coincidence, and we will now connect the above result of our number-theoretic detour with lattice triangles and $2$-collinear integers.

\begin{theorem}Only $1$, $2$, $4$, and $7$ are $2$-collinear integers.\end{theorem}
\begin{proof}Let us define $D_{2k+1}$ to be the subset of $\{0,1,2,\dots,2k\}$ such that $a\in D_{2k+1}$ iff $\gcd(a,2k+1)=\gcd(a-1,2k+1)=1.$ Following \cite{reznick}, we obtain that $$|D_{2k+1}|=\phi(2k+1,1).$$

Now suppose that $a\in D_{2k+1},$ so that $\gcd(a,2k+1)=\gcd(a-1,2k+1)=1.$ By the Euclidean Algorithm, it follows that $\gcd(2k+1-a,2k+1)=\gcd(-a,2k+1)=1$ and $\gcd(2k+2-a,2k+1)=\gcd(-a+1,2k+1)=1,$ so $2k+2-a\in D_{2k+1}$ as well; to avoid bounding issues, note that both $0$ and $1$ are not in $D_{2k+1}$, so $2k+2-a\in\{0,1,2,\dots,2k\}$.

The above implies that if $a\in D_{2k+1}$, then $2k+2-a\in D_{2k+1}$, implying a bijection between $$\{1,2,3,\dots,k\}\cap D_{2k+1}$$ and $$\{k+2,k+3,k+4,\dots,2k+1\}\cap D_{2k+1}.$$ The Euclidean Algorithm again implies that $\gcd(k+1,2k+1)=\gcd(k,2k+1)=1,$ so it follows that $k+1\in D_{2k+1}$. Thus $$|\{1,2,3,\dots,k\}\cap D_{2k+1}|+|\{k+2,k+3,k+4,\dots,2k+1\}\cap D_{2k+1}|=|D_{2k+1}|-1,$$ and the bijective property implies that $$|\{1,2,3,\dots,k\}\cap D_{2k+1}|=|\{k+2,k+3,k+4,\dots,2k+1\}\cap D_{2k+1}|=\frac{|D_{2k+1}|-1}{2}.$$

With the above, we focus on the set $\{3,4,5,\dots,k\}\cap D_{2k+1}.$  For ease of notation, let $$S_k=\{3,4,5,\dots,k\}.$$  As established earlier, note that $1\notin D_{2k+1}$, while $$\gcd(1,2k+1)=\gcd(2,2k+1)=1$$ and $\gcd(0,2k+1)=2k+1$ gives $\{1,2\}\cap D_{2k+1}=\{2\}$, and therefore 
\begin{align*}
|S_k \cap D_{2k+1}|=|\{1,2,3,\dots,k\}\cap D_{2k+1}|-1 =\frac{|D_{2k+1}|-3}2.
\end{align*}

Now suppose that $$k\notin\{1,2,4,7\},$$ so that $2k+1$ is odd and is not a member of the set $\{3,5,9,15\}.$ By Lemma~\ref{stg5}, it follows that $|D_{2k+1}|=\phi(2k+1,1)\ge 5,$ so $$|S_k \cap D_{2k+1}|=\frac{|D_{2k+1}|-3}2\ge\frac{5-3}{2}=1.$$
Hence there must exist some $a\in S_k$ with $\gcd(a,2k+1)=\gcd(a-1,2k+1)=1,$ so it follows by Theorem~\ref{badk} that $k$ is not $2$-collinear. Combined with Corollary~\ref{1247}, we obtain that the set of all $2$-collinear integers are $1$, $2$, $4$, and $7$.
\end{proof}

\section{Conclusion}

In this paper, we have demonstrated that the only integers $k$ such that a lattice triangle $T$ satisfying $$B(T)=3\text{ and }I(P)=k$$ has all $k$ interior points necessarily collinear are $$k=1,~2,~4,\text{ and }7.$$  In future research, we plan to consider extending this work to more general lattice $n$-gons, and to lattice polytopes in higher dimensions.
\makeatletter
\renewcommand{\@biblabel}[1]{[#1]\hfill}
\makeatother

\end{document}